\documentclass[12pt]{amsart}
\usepackage{amssymb}
\usepackage[all]{xy}
\usepackage{tikz}
\usepackage{amscd}
\usepackage{amsmath}
\usepackage{mathrsfs,graphicx,amsthm, mathtools}
\numberwithin{equation}{section}
\input xy
\xyoption{all}

\newtheorem{lemma}{Lemma}[section]

\newtheorem{theorem}[lemma]{Theorem}
\newtheorem{proposition}[lemma]{Proposition}

\theoremstyle{definition}
\newtheorem{remark}[lemma]{Remark}
\newtheorem{definition}[lemma]{Definition}

\DeclareMathOperator{\Mod}{Mod}
\DeclareMathOperator{\modd}{mod}

\DeclareMathOperator{\Add}{Add}

\DeclareMathOperator{\Hom}{Hom}

\DeclareMathOperator{\Ext}{Ext}

\DeclareMathOperator{\Lim}{Lim}

\DeclareMathOperator{\fp}{fp}
\DeclareMathOperator{\Flatt}{Flat}
\DeclareMathOperator{\ind}{ind}
\DeclareMathOperator{\Ind}{Ind}
\DeclareMathOperator{\Eff}{Eff}
\DeclareMathOperator{\End}{End}
\DeclareMathOperator{\eff}{eff}

\DeclareMathOperator{\A}{A}

\newtheorem{question}[lemma]{Question}

\newcounter{diagram}
\numberwithin{diagram}{section}

\newtheorem*{theorem a*}{Theorem A}
\newtheorem*{theorem b*}{Theorem B}
\newtheorem*{theorem c*}{Theorem C}
\newtheorem*{corollary a*}{Corollary A}
\newtheorem*{question a*}{Question A}
\begin{document}
	
	\title{The completion of $d$-abelian categories}
	
	\author{Ramin Ebrahimi}
	\address{Department of Pure Mathematics\\
		Faculty of Mathematics and Statistics\\
		University of Isfahan\\
		P.O. Box: 81746-73441, Isfahan, Iran\\ and School of Mathematics, Institute for Research in Fundamental Sciences (IPM), P.O. Box: 19395-5746, Tehran, Iran}
	\email{r.ebrahimi@ipm.ir / ramin.ebrahimi1369@gmail.com}
	
	\author{Alireza Nasr-Isfahani}
	\address{Department of Pure Mathematics\\
		Faculty of Mathematics and Statistics\\
		University of Isfahan\\
		P.O. Box: 81746-73441, Isfahan, Iran\\ and School of Mathematics, Institute for Research in Fundamental Sciences (IPM), P.O. Box: 19395-5746, Tehran, Iran}
	\email{nasr$_{-}$a@sci.ui.ac.ir / nasr@ipm.ir}

	

		\begin{abstract}
		Let $A$ be a finite-dimensional algebra, and $\mathfrak{M}$ be a $d$-cluster tilting subcategory of $\modd A$. From the viewpoint of higher homological algebra, a natural question to ask is when $\mathfrak{M}$ induces a $d$-cluster tilting subcategory in $\Mod A$. In this paper, we investigate this question in a more general form. We consider $\mathcal{M}$ as a small $d$-abelian category, known to be equivalent to a $d$-cluster tilting subcategory of an abelian category $\mathcal{A}$. The completion of $\mathcal{M}$, denoted by $\Ind(\mathcal{M})$, is defined as the universal completion of $\mathcal{M}$ with respect to filtered colimits. We explore $\Ind(\mathcal{M})$ and demonstrate its equivalence to the full subcategory $\mathcal{L}_d(\mathcal{M})$ of $\Mod \mathcal{M}$, comprising left $d$-exact functors.
		Notably, while $\Ind(\mathcal{M})$ as a subcategory of $\frac{\Mod \mathcal{M}}{\Eff(\mathcal{M})}$, satisfies all properties of a $d$-cluster tilting subcategory except $d$-rigidity, it falls short of being a $d$-cluster tilting category. For a $d$-cluster tilting subcategory $\mathfrak{M}$ of $\modd A$, $\overrightarrow{\mathfrak{M}}$, consists of all filtered colimits of objects from $\mathfrak{M}$, is a generating-cogenerating, functorially finite subcategory of $\Mod A$. The question of whether $\mathfrak{M}$ is a $d$-rigid subcategory remains unanswered. However, if it is indeed $d$-rigid, it qualifies as a $d$-cluster tilting subcategory. In the case $d=2$, employing cotorsion theory, we establish that $\overrightarrow{\mathfrak{M}}$ is a $2$-cluster tilting subcategory if and only if $\mathfrak{M}$ is of finite type. Thus, the question regarding whether $\overrightarrow{\mathfrak{M}}$ is a $d$-cluster tilting subcategory of $\Mod A$ appears to be equivalent to the Iyama's qestion about the finiteness of $\mathfrak{M}$. Furthermore, for general $d$, we address the problem and present several equivalent conditions for the Iyama's question.
	\end{abstract}
	
	\maketitle

	%%%%%%%%%%%%%%%%%%%%%%%%%%%%%%%%%%%%%%

	\section{Introduction}
	For a positive integer $d$, a $d$-cluster tilting subcategory of abelian and exact categories was defined by Iyama to construct higher-dimensional analogs of Auslander-Reiten theory and the Auslander correspondence.
	$d$-abelian categories were introduced by Jasso to axiomatize $d$-cluster tilting subcategories of abelian categories \cite{J}.
	Jasso proved that every $d$-cluster tilting subcategory of an abelian category is a $d$-abelian category. Recently, it has been shown that any $d$-abelian category is
	equivalent to a $d$-cluster tilting subcategory of an abelian category \cite{EN1,Kv}.
	When $d\geq 2$, the study of these categories and $d$-abelian categories is sometimes called "Higher Homological Algebra".
	
	For an additive category $\mathcal{M}$, we denote by $\Ind(\mathcal{M})$ the completion (ind-completion or indization in \cite{KS}) of $\mathcal{M}$. $\Ind(\mathcal{M})$ is defined as the completion of $\mathcal{M}$ with respect to filtered colimits. Indeed, $\Ind(\mathcal{M})$ is defined by the following universal property: $\Ind(\mathcal{M})$ is an additive category with filtered colimits with an additive functor
	\begin{equation}
		\mathbb{I}:\mathcal{M}\longrightarrow \Ind(\mathcal{M}) \notag
	\end{equation}
	such that any other additive functor $\mathcal{M}\rightarrow \mathcal{C}$, where $\mathcal{C}$ is an additive category with filtered colimits, must factor uniquely trough $\mathbb{I}$. The completion of every small category exists and is unique up to equivalence. One way to construct $\Ind(\mathcal{M})$ is as follows: 
	Let $\mathcal{M}$ be a small additive category. We denote by $\Mod\mathcal{M}$ the category of all additive
	contravariant functors from $\mathcal{M}$ to the category of all abelian groups. The Yoneda functor
	\begin{equation}
		\mathcal{Y}:\mathcal{M}\longrightarrow \rm Mod\mathcal{M} \notag
	\end{equation}
	is the free completion of $\mathcal{M}$ (complete with respect to all colimits). Let $\Ind(\mathcal{M})$ be the full subcategory of $\Mod \mathcal{M}$ whose objects are filtered
	colimits of representable functors. Hence, the Yoneda functor induces a functor
	\begin{equation}
		\mathbb{I}:\mathcal{M}\longrightarrow \Ind(\mathcal{M}) \notag
	\end{equation}
	that is the completion of $\mathcal{M}$.
	
	Let $\mathcal{M}$ be a small $d$-abelian category, and $\Eff(\mathcal{M})$ be the category of all  weakly effaceable functors. By \cite{E}, we have the following localization of abelian categories:
	\begin{center}
		\begin{tikzpicture}
			\node (X1) at (-3,0) {$\Eff(\mathcal{M})$};
			\node (X2) at (0,0) {$\Mod\mathcal{M}$};
			\node (X3) at (3,0) {$\frac{\Mod\mathcal{M}}{\Eff(\mathcal{M})}$};
			\node (p) at (-1.5,-0.7) {$p$};
			\node (r) at (1.5,-0.7) {$r$};
			\draw [->,thick] (X1) -- (X2) node [midway,above] {$i$};
			\draw [->,thick] (X2) -- (X3) node [midway,above] {$e$};
			\draw [->,thick] (X2) to [out=225,in=315] (X1) node [midway,left] {};
			\draw [->,thick] (X3) to [out=225,in=315] (X2) node [midway,left] {};
		\end{tikzpicture}
	\end{center}
	Let $\mathcal{L}_d(\mathcal{M})$ be the full subcategory of $\rm Mod\mathcal{M}$ consisting of left $d$-exact functors. We show that $\Ind(\mathcal{M})\simeq\mathcal{L}_d(\mathcal{M})$, see Proposition 3.3. First, for an arbitrary small $d$-abelian category $\mathcal{M}$, we prove some basic features of $\Ind(\mathcal{M})$.\\

	\textbf{Theorem A.} (Theorem \ref{3.4})
		Let $\mathcal{M}$ be a small $d$-abelian category. By the above notations we have
		\begin{itemize}
			\item[(i)]
			$\Ind(\mathcal{M})$ is a generating-cogenerating functorialy finite subcategory of an abelian category $\frac{\Mod\mathcal{M}}{\Eff(\mathcal{M})}$.
			\item[(ii)]
			For every $F\in \Mod\mathcal{M}$, there exists
			an exact sequence
			\begin{equation}
				0\rightarrow L_d\rightarrow L_{d-1}\rightarrow \cdots\rightarrow L_1\rightarrow F\rightarrow 0 \notag
			\end{equation}
			in $\mathcal{L}_1(\mathcal{M})$ such that $L_i$'s are in $\mathcal{L}_d(\mathcal{M})$.
			\item[(iii)]
			We have $\Ind(\mathcal{M})^{_1\bot_{d-1}}\subseteq \Ind(\mathcal{M})$.\\
		\end{itemize}
	
	Note that if we know $\Ind(\mathcal{M})$ is a $d$-rigid subcategory of
	$\frac{\Mod\mathcal{M}}{\Eff(\mathcal{M})}$, then part (ii) of the above theorem implies that $^{_1\bot_{d-1}}\Ind(\mathcal{M})\subseteq \Ind(\mathcal{M})$.
	Thus $\Ind(\mathcal{M})$ is a $d$-cluster tilting subcategory of $\frac{\Mod\mathcal{M}}{\Eff(\mathcal{M})}$ if and only if it is a $d$-rigid subcategory. By using a result from \cite{E2}, we will show that $\Ind(\mathcal{M})=\mathcal{L}_d(\mathcal{M})$ is $d$-rigid in $\frac{\Mod\mathcal{M}}{\Eff(\mathcal{M})}$ if and only if it is $d$-rigid in $\Mod\mathcal{M}$.
	We don't know yet if this condition is satisfied or not in general. But since $\Ind(\mathcal{M})$ is equal to $\Flatt(\mathcal{M})$, the subcategory of flat functors in $\Mod\mathcal{M}$, we have the following corollary.\\

	\textbf{Corollary A.} 
		Let $\mathcal{M}$ be a small $d$-abelian category, that is a perfect category (i.e. flat functors and projective functors coincide.). Then $\Ind(\mathcal{M})$ is a $d$-cluster tilting subcategory of $\frac{\Mod\mathcal{M}}{\Eff(\mathcal{M})}$.\\
	
	Motivated by the above results, we conjecture that the following question has a positive answer.\\
	
	\textbf{Question A.} (Question \ref{3.5})
		Let $\mathcal{M}$ be a small $d$-abelian category. Then $\Ind(\mathcal{M})$ is a $d$-rigid subcategory of $\Mod\mathcal{M}$.\\
	
	Let $A$ be a finite-dimensional algebra and $\mathfrak{M}$ be a $d$-cluster tilting subcategory of $\modd A$. By the above results, we have that $\overrightarrow{\mathfrak{M}}$, which consists of all filtered colimits of objects from $\mathfrak{M}$, is a generating-cogenerating and functorially finite subcategory of $\Mod A$. If it is $d$-rigid then it is indeed a $d$-cluster tilting subcategory.
	For the case $d=2$, by using cotorsion theory, we prove the following result.\\

	\textbf{Theorem B.} (Theorem \ref{4.7})
		Let $A$ be a finite-dimensional algebra and $\mathfrak{M}$ be a $2$-cluster tilting subcategory of $\modd A$. The following are equivalent.
		\begin{itemize}
			\item[(a)] $\overrightarrow{\mathfrak{M}}$ is a $2$-cluster tilting subcategory of $\Mod A$.
			\item[(b)] $\overrightarrow{\mathfrak{M}}$ is a $2$-rigid subcategory of $\Mod A$.
			\item[(c)] $\mathfrak{M}$ has finitely many indecomposables, up to isomorphism.\\
		\end{itemize}

The question of the finiteness of d-cluster tilting subcategories for $d\geq2$, which was one of the first questions posed by Iyama \cite{I3}, is still open. To date, all known d-cluster tilting subcategories with $d\geq2$ are of finite type. Therefore, the question of whether $\overrightarrow{\mathfrak{M}}$ is a $d$-cluster tilting subcategory of $\Mod A$ seems to be equivalent to Iyama's question about the finiteness of $\mathfrak{M}$.
	
	This paper is organized as follows. In section 2, we recall the definitions of $d$-abelian categories, $d$-cluster tilting subcategories, and some of their basic properties. Also, for a $d$-abelian category $\mathcal{M}$, we recall some results about the functor categories $\Mod\mathcal{M}$ and $\modd\mathcal{M}$ and their localizations $\frac{\Mod\mathcal{M}}{\Eff(\mathcal{M})}$ and $\frac{\modd\mathcal{M}}{\eff(\mathcal{M})}$, with respect to weakly effaceable and effaceable functors, respectively.
	In section 3, we prove that the completion of a small $d$-abelian category $\mathcal{M}$ is equivalent to the full subcategory $\mathcal{L}_d(\mathcal{M})$ of left $d$-exact functors in $\Mod\mathcal{M}$ and study properties of this category. Also, we discuss the situation when $\Ind(\mathcal{M})$ is a $d$-abelian category. In section 4, for a finite-dimensional algebra $A$ we consider the case of $d$-cluster tilting subcategories of $\modd A$.
	
	\subsection{Notation}
	Throughout this paper, all categories and functors are additive. $d$ always denotes a fixed positive integer and $\mathcal{M}$ is a fixed small $d$-abelian category.

	%%%%%%%%%%%%%%%%%%%%%%%%%%%%%%%%%%%%%%
	
	\section{preliminaries}
	In this section, we recall the definitions of $d$-abelian categories, $d$-cluster tilting subcategories, and some of their basic properties. Also, for a $d$-abelian category $\mathcal{M}$, we recall some results about the functor categories $\Mod\mathcal{M}$ and $\modd\mathcal{M}$.
	
	\subsection{$d$-abelian categories}
	Let $\mathcal{M}$ be an additive category and $f:A\rightarrow B$ a morphism in $\mathcal{M}$. A {\em weak cokernel} of $f$ is a morphism $g:B\rightarrow C$ such that for all $C^{\prime} \in \mathcal{M}$  the sequence of abelian groups
	\begin{equation}
		\Hom_{\mathcal{M}}(C,C')\overset{(g,C')}{\longrightarrow} \Hom_{\mathcal{M}}(B,C')\overset{(f,C')}{\longrightarrow} \Hom_{\mathcal{M}}(A,C') \notag
	\end{equation}
	is exact. The concept of {\em weak kernel} is defined dually.
	
	Let $\alpha^0:X^0 \rightarrow X^1$ be a morphism in $\mathcal{M}$. A {\em $d$-cokernel} of $\alpha^0$ is a sequence
	\begin{equation}
		(\alpha^1, \ldots, \alpha^d): X^1 \overset{\alpha^1}{\longrightarrow} X^2 \overset{\alpha^2}{\longrightarrow}\cdots \overset{\alpha^{d-1}}{\longrightarrow} X^d \overset{\alpha^d}{\longrightarrow} X^{d+1} \notag
	\end{equation}
	of objects and morphisms in $\mathcal{M}$ such that, for each $Y\in \mathcal{M}$
	the induced sequence of abelian groups
	\begin{align}
		0 \rightarrow \Hom_{\mathcal{M}}(X^{d+1},Y) \rightarrow \Hom_{\mathcal{M}}(X^d,Y) \rightarrow\cdots\rightarrow \Hom_{\mathcal{M}}(X^1,Y) \rightarrow \Hom_{\mathcal{M}}(X^0,Y) \notag
	\end{align}
	is exact. Equivalently, the sequence $(\alpha^1, \ldots, \alpha^d)$ is a $d$-cokernel of $\alpha^0$ if for all $1\leq k\leq d-1$
	the morphism $\alpha^k$ is a weak cokernel of $\alpha^{k-1}$, and $\alpha^d$ is moreover a cokernel of $\alpha^{d-1}$. The concept of {\em $d$-kernel} of a morphism is defined dually \cite[Definition 2.2]{J}.
	\begin{definition}(\cite[Definition 2.4]{L})\label{d1}
		Let $\mathcal{M}$ be an additive category. A {\em left $d$-exact sequence} in $\mathcal{M}$ is a complex
		\begin{equation}
			X^0\overset{\alpha^0}{\longrightarrow}X^1 \overset{\alpha^1}{\longrightarrow} X^2 \overset{\alpha^2}{\longrightarrow}\cdots \overset{\alpha^{d-1}}{\longrightarrow} X^d \overset{\alpha^d}{\longrightarrow} X^{d+1} \notag
		\end{equation}
		such that $(\alpha^0, \ldots, \alpha^{d-1})$ is a $d$-kernel of $\alpha^d$. The concept of {\em right $d$-exact sequence} is defined dually. A {\em $d$-exact sequence} is a sequence which is both right $d$-exact and left $d$-exact.
	\end{definition}
	
	Let $\mathcal{M}$ be a category and $A$ be an object of $\mathcal{M}$. A morphism $e\in \Hom_{\mathcal{M}}(A, A)$ is called {\em idempotent} if $e^2 = e$. $\mathcal{M}$ is called {\em idempotent complete}
	if for every idempotent $e\in \Hom_{\mathcal{M}}(A, A)$ there exist an object $B$ and morphisms $f\in\Hom_{\mathcal{M}}(A, B)$ and
	$g\in\Hom_{\mathcal{M}}(B, A)$ such that $gf = e$ and $fg = 1_B$ \cite[Page 61]{Fr}.
	
	\begin{definition}$($\cite[Definition 3.1]{J}$)$\label{2.2}
		A {\em $d$-abelian} category is an additive category $\mathcal{M}$ which satisfies the following axioms.
		\begin{itemize}
			\item[$(\A 0)$]
			The category $\mathcal{M}$ is idempotent complete.
			\item[$(\A 1)$]
			Every morphism in $\mathcal{M}$ has a $d$-kernel and a $d$-cokernel.
			\item[$(\A 2)$]
			For every monomorphism $\alpha^0:X^0 \rightarrow X^1$ in $\mathcal{M}$ and for every $d$-cokernel $(\alpha^1, \ldots, \alpha^d)$ of $\alpha^0$, the following sequence is $d$-exact:
			\begin{equation}
				X^0\overset{\alpha^0}{\longrightarrow}X^1 \overset{\alpha^1}{\longrightarrow} X^2 \overset{\alpha^2}{\longrightarrow}\cdots \overset{\alpha^{d-1}}{\longrightarrow} X^d \overset{\alpha^d}{\longrightarrow} X^{d+1}. \notag
			\end{equation}
			\item[$(\A 2^{op})$]
			For every epimorphism $\alpha^d:X^d \rightarrow X^{d+1}$ in $\mathcal{M}$ and for every $d$-kernel $(\alpha^0, \ldots, \alpha^{d-1})$ of $\alpha^d$, the following sequence is $d$-exact:
			\begin{equation}
				X^0\overset{\alpha^0}{\longrightarrow}X^1 \overset{\alpha^1}{\longrightarrow} X^2 \overset{\alpha^2}{\longrightarrow}\cdots \overset{\alpha^{d-1}}{\longrightarrow} X^d \overset{\alpha^d}{\longrightarrow} X^{d+1}.
			\end{equation}
		\end{itemize}
	\end{definition}
	
	A subcategory $\mathcal{B}$ of an abelian category $\mathcal{A}$ is called {\em cogenerating} if for every object
	$X\in \mathcal{A}$ there exists an object $Y\in\mathcal{B}$ and a monomorphism $X\rightarrow Y$. The concept of
	{\em generating} subcategory is defined dually \cite[Page 69]{Fr}.
	
	Let $\mathcal{A}$ be an additive category and $\mathcal{B}$ be a full subcategory of $\mathcal{A}$. $\mathcal{B}$ is called
	{\em covariantly finite in} $\mathcal{A}$ if for every $A\in \mathcal{A}$ there exist an object $B\in\mathcal{B}$ and a morphism
	$f : A\rightarrow B$ such that, for all $B'\in\mathcal{B}$, the sequence of abelian groups $\Hom_{\mathcal{A}}(B, B')\rightarrow \Hom_{\mathcal{A}}(A, B')\rightarrow 0$ is exact. Such a morphism $f$ is called a {\em left $\mathcal{B}$-approximation of} $A$. The notions of {\em contravariantly finite subcategory} of $\mathcal{A}$ and {\em right $\mathcal{B}$-approximation} are defined dually. A {\em functorially
		finite subcategory} of $\mathcal{A}$ is a subcategory which is both covariantly and contravariantly finite
	in $\mathcal{A}$ \cite[Page 113]{AR}.
	
	\begin{definition}$($\cite[Definition 3.14]{J}$)$
		Let $\mathcal{A}$ be an abelian category and $\mathcal{M}$ be a generating-cogenerating full subcategory of $\mathcal{A}$. $\mathcal{M}$ is called a {\em $d$-cluster tilting subcategory} of $\mathcal{A}$ if $\mathcal{M}$ is functorially finite in $\mathcal{A}$ and
		\begin{align}
			\mathcal{M}& = \{ X\in \mathcal{A} \mid \forall i\in \{1, \ldots, d-1 \}, \Ext_{\mathcal{A}}^i(X,\mathcal{M})=0 \}\notag \\
			& =\{ X\in \mathcal{A} \mid \forall i\in \{1, \ldots, d-1 \}, \Ext_{\mathcal{A}}^i(\mathcal{M},X)=0 \}.\notag
		\end{align}
		
		Note that $\mathcal{A}$ itself is the unique $1$-cluster tilting subcategory of $\mathcal{A}$.
	\end{definition}
	
	\begin{remark}$($\cite[Remark 3.15]{J}$)$\label{2.5}
		Let $\mathcal{A}$ be an abelian category and $\mathcal{M}$ be a $d$-cluster tilting subcategory of $\mathcal{A}$. Since $\mathcal{M}$ is a generating-cogenerating subcategory of $\mathcal{A}$, for each $A\in \mathcal{A}$, every left $\mathcal{M}$-approximation of $A$ is a monomorphism and every right $\mathcal{M}$-approximation of $A$ is an epimorphism.
	\end{remark}
	
	A full subcategory $\mathcal{M}$ of an abelian category $\mathcal{A}$ is called {\em $d$-rigid}, if for every two objects $M,N\in \mathcal{M}$ and for every $k\in \{1,\ldots,d-1\}$, we have $\Ext_{\mathcal{A}}^k(M,N)=0$ \cite[Page 443]{Bel}. Any $d$-cluster tilting subcategory $\mathcal{M}$ of an abelian category $\mathcal{A}$ is $d$-rigid.
	
	\begin{theorem}\label{2.6}
		\begin{itemize}
			\item[(i)]
			Let $\mathcal{A}$ be an abelian category and $\mathcal{M}$ be a $d$-cluster tilting subcategory of $\mathcal{A}$. Then $\mathcal{M}$ is a $d$-abelian category \cite[Theorem 3.16]{J}.
			\item[(ii)]
			Conversely, for every small $d$-abelian category $\mathcal{M}$ there exists an abelian category $\mathcal{A}$, such that $\mathcal{M}$ is equivalent to a $d$-cluster tilting subcategory of $\mathcal{A}$ \cite[Theorem 4.7]{EN1} and \cite[Corollary 1.3]{Kv}.
		\end{itemize}
	\end{theorem}
	
	Let $\mathcal{M}$ be a small $d$-abelian category. Recall that $\Mod\mathcal{M}$ is
	the category of all additive contravariant functors from $\mathcal{M}$ to the category of all abelian
	groups. It is an abelian category with all limits and colimits, which are defined point-wise. A functor $F\in \Mod\mathcal{M}$ is called {\em finitely presented} or {\em coherent}, if there exists an exact sequence
	\begin{center}
		$\Hom_{\mathcal{M}}(-,X)\rightarrow \Hom_{\mathcal{M}}(-,Y)\rightarrow F\rightarrow 0$
	\end{center}
	in $\Mod\mathcal{M}$. We denote by $\modd\mathcal{M}$ the full subcategory of $\Mod\mathcal{M}$ consist of all finitely presented functors. It is a well known result that $\modd\mathcal{M}$ is an abelian category if and only if $\mathcal{M}$ has weak kernels and in this case the inclusion functor $\modd\mathcal{M}\hookrightarrow \Mod\mathcal{M}$ is an exact functor \cite[Theorem 1.4]{Fr2}.
	By Yoneda's lemma, representable functors are projective and the direct sum of all representable functors $\bigoplus_{X\in \mathcal{M}}\Hom_{\mathcal{M}}(-,X)$ is a generator for $\Mod\mathcal{M}$. Thus $\Mod\mathcal{M}$ is a Grothendieck category \cite[Proposition 5.21]{Fr}.
	
	\begin{definition}$($\cite[Definition 3.2]{EN1}$)$ \label{2.7}
		Let $\mathcal{M}$ be a $d$-abelian category, $\mathcal{A}$ an abelian category and $F:\mathcal{M}\rightarrow \mathcal{A}$ a contravariant functor. We say that $F$ is a
		\begin{itemize}
			\item[(i)]
			{\em Left $d$-exact functor} if for any left $d$-exact sequence $X^0 \rightarrow X^1 \rightarrow \cdots \rightarrow X^d \rightarrow X^{d+1} $ in $\mathcal{M}$, $0 \rightarrow F(X^{d+1}) \rightarrow F(X^d) \rightarrow \cdots \rightarrow F(X^1) \rightarrow F(X^0)$ is an exact sequence of $\mathcal{A}$.
			\item[(ii)]
			{\em Right $d$-exact functor} if for any right $d$-exact sequence $X^0 \rightarrow X^1 \rightarrow \cdots \rightarrow X^d \rightarrow X^{d+1} $ in $\mathcal{M}$, $F(X^{d+1}) \rightarrow F(X^d) \rightarrow \cdots \rightarrow F(X^1) \rightarrow F(X^0)\rightarrow 0$ is an exact sequence of $\mathcal{A}$.
			\item[(iii)]
			{\em $d$-exact functor} if for any $d$-exact sequence $X^0 \rightarrow X^1 \rightarrow \cdots \rightarrow X^d \rightarrow X^{d+1} $ in $\mathcal{M}$, $0 \rightarrow F(X^{d+1}) \rightarrow F(X^d) \rightarrow \cdots \rightarrow F(X^1) \rightarrow F(X^0)\rightarrow 0$ is an exact sequence of $\mathcal{A}$.
		\end{itemize}
		The covariant left $d$-exact (resp., right $d$-exact, $d$-exact) functors are defined similarly.
	\end{definition}
	
	\begin{definition}\label{2.8}
		\begin{itemize}
			\item[(i)]
			A functor $F\in \Mod\mathcal{M}$ is called {\em weakly effaceable}, if for each object $X\in \mathcal{M}$ and $x\in F(X)$
			there exists an epimorphism $f : Y \rightarrow X$ such that $F(f)(x) = 0$. We denote by
			$\Eff(\mathcal{M})$ the full subcategory of all weakly effaceable functors.
			\item[(ii)]
			A functor $F\in \modd\mathcal{M}$ is called {\em effaceable}, if there exists an exact sequence
			\begin{center}
				$\Hom_{\mathcal{M}}(-,Y)\rightarrow \Hom_{\mathcal{M}}(-,X)\rightarrow F\rightarrow 0$
			\end{center}
			such that $Y \rightarrow X$ is an epimorphism. We denote by
			$\eff(\mathcal{M})$ the full subcategory of all effaceable functors.
		\end{itemize}
	\end{definition}
	
	\subsection{Localisation with respect to effaceable functors}
	In this section we consider the localization of abelian categories $\Mod \mathcal{M}$ and $\modd \mathcal{M}$ with respect to so called the subcategories of weakly effaceable and effaceable functors, respectively.
	For more information about the localisation theory of abelian categories we refer the reader to Gabriel thesis \cite{G} or the recent book of Krause \cite{KrB}.
	
Let $\mathcal{A}$ be an abelian category and $\mathcal{X}$ be a subcategory of $\mathcal{A}$. For integers $0\leq i\leq k$ we denote by $\mathcal{X}^{_i{\bot}_k}$ the full subcategory of $\mathcal{A}$ defined by $\mathcal{X}^{_i{\bot}_k} = \{A \in \mathcal{A} | \Ext_{\mathcal{A}}^{i,...,k}(\mathcal{X}, A) = 0\}$. $^{_i{\bot}_k}\mathcal{X}$ is defined similarly \cite[Page 72]{Ps}.
Note that for a localising subcategory $\mathcal{C}$ of $\mathcal{A}$ we have $\mathcal{C}^{_0{\bot}_1}=\mathcal{C}^{\bot}$, the subcategory of $\mathcal{C}$-closed objects \cite[Lemma 2.2.4]{KrB}.
	
	The following proposition from \cite{E2} can be seen as a variation of \cite[Proposition 3.4 and Theorem 3.10]{Ps}.
	
	\begin{proposition}$($\cite[Proposition 2.3]{E2}$)$\label{2.10}
		Let $\mathcal{A}$ be a Grothendieck category and $\mathcal{C}$ be a localising subcategory of $\mathcal{A}$. For an object $A\in \mathcal{A}$ and a non-negative integer $k$ the following statements are equivalent.
		\begin{itemize}
			\item[$(i)$]
			$A\in \mathcal{C}^{{\bot}_{k+1}}$.
			\item[$(ii)$]
			There exists an injective coresolution
			\begin{equation}\label{ri}
				0\rightarrow A\rightarrow r(I^0) \rightarrow r(I^1)\rightarrow \cdots\rightarrow r(I^{k+1})
			\end{equation}
			for $A$, where $I^i$'s are injective objects in ${\mathcal{A}}/{\mathcal{C}}$.
			\item[$(iii)$]
			The natural map $e_{X,A}^i:\Ext_{\mathcal{A}}^i(X,A)\rightarrow\Ext_{{\mathcal{A}}/{\mathcal{C}}}^i(e(X),e(A))$ is invertible, for every $X\in \mathcal{A}$ and every $i\in\{0,1,\ldots,k\}$.
		\end{itemize}
	\end{proposition}
	
	\begin{proposition}\label{2.11}
		Let $\mathcal{M}$ be a small $d$-abelian category. Then the following statements hold.
		\begin{itemize}
			\item[(i)]
			$\Eff(\mathcal{M})$ is a localising subcategory of $\Mod\mathcal{M}$.
			\item[(ii)]
			For every $k\in \{1,\ldots, d\}$, $\Eff(\mathcal{M})^{_0\bot_k}=\mathcal{L}_k(\mathcal{M})$. In particular
			\begin{align*}
				\mathcal{L}_1(\mathcal{M})&=\Eff(\mathcal{M})^{\bot}\simeq\dfrac{\Mod\mathcal{M}}{\Eff(\mathcal{M})},\\
				\mathcal{L}_d(\mathcal{M})&=\Eff(\mathcal{M})^{_0\bot_d}.
			\end{align*}
			\item[(iii)]
			$\eff(\mathcal{M})$ is a localising subcategory of $\modd\mathcal{M}$.
			\item[(iv)]
			The localisation functor $\mathcal{M}\longrightarrow \frac{\modd\mathcal{M}}{\eff(\mathcal{M})}$ is fully faithful and the essential image of this functor is a $d$-cluster tilting subcategory.
		\end{itemize}
		\begin{proof}
			(i) and (ii) was proved in \cite[Propositions 3.5 and 3.6]{E}. (iii) and (iv) was proved in \cite[Proposition 4.4 and Theorm 4.7]{EN1}.
		\end{proof}
	\end{proposition}
	
	By the above proposition $\mathcal{L}_1(\mathcal{M})$ is a Grothendieck category. Thus $\mathcal{L}_1(\mathcal{M})$ has injective envelopes. We will need the following easy lemma.
	
	\begin{lemma}\label{2.12}
		Every injective object in $\mathcal{L}_1(\mathcal{M})$ is a $d$-exact functor.
		\begin{proof}
			Every injective object $E\in \mathcal{L}_1(\mathcal{M})$ is also an injective object in $\Mod\mathcal{M}$. So by \cite[Proposition 3.5]{EN3} it is a right $d$-exact functor. Since $E\in \mathcal{L}_1(\mathcal{M})$, the result follows.
		\end{proof}
	\end{lemma}
	
	The following lemma is the combination of Propositions 2.10 and 2.11.
	\begin{lemma}\label{2.13}
		Let $\mathcal{M}$ be a small $d$-abelian category and $k\in \{1,\ldots,d\}$ be a positive integer. The following statements are equivalent for an object $F\in \Mod\mathcal{M}$.
		\begin{itemize}
			\item[(i)]
			$F\in \mathcal{L}_k(\mathcal{M})$.
			\item[(ii)]
			The natural map
			$e_{G,F}^i:\Ext_{\Mod\mathcal{M}}^i(G,F)\rightarrow\Ext_{\Mod\mathcal{M}/{\Eff(\mathcal{M})}}^i(G,F)$
			is invertible, for every $G\in \Mod\mathcal{M}$ and every $i\in\{0,1,\ldots,k\}$.
		\end{itemize}
	\end{lemma}
	
	We denote by $H:\mathcal{M}\rightarrow\mathcal{L}_1(\mathcal{M})$ the composition of the Yoneda functor $\mathcal{M}\rightarrow \Mod\mathcal{M}$ and the localisation functor $\Mod\mathcal{M}\rightarrow \frac{\Mod\mathcal{M}}{\Eff(\mathcal{M})}\simeq \Eff(\mathcal{M})^{\bot}=\mathcal{L}_1(\mathcal{M})$.
	Thus $H(X)=(-,X):\mathcal{M}^{op}\rightarrow \rm Ab$.
	
	\subsection{$\Ind$-completion and locally finitely presented categories}
	Let $\mathcal{B}$ be an additive category with all direct limits.
	\begin{itemize}
		\item[(i)]
		An object $X\in \mathcal{B}$ is called
		{\em finitely presented (finitely generated)} provided that for every direct limit, $\underrightarrow{\Lim}_{i\in \mathcal{I}}X_i$ in
		$\mathcal{B}$ the natural morphism
		\begin{center}
			$\underrightarrow{\Lim}\Hom_{\mathcal{B}}(X,Y_i)\rightarrow \Hom_{\mathcal{B}}(X,\underrightarrow{\Lim}Y_i)$
		\end{center}
		is an isomorphism (a monomorphism).
		The full subcategory of finitely presented objects of $\mathcal{A}$ is denoted by $\fp(\mathcal{A})$.
		\item[(ii)]
		$\mathcal{B}$ is called a {\em locally finitely presented category} if
		$\fp(\mathcal{B})$ is skeletally small and every object in $\mathcal{B}$ is a direct limit of objects in $\fp(\mathcal{B})$.
		\item[(iii)]
		Let $\mathcal{B}$ be a locally finitely presented abelian category. $\mathcal{B}$ is said to be {\em locally coherent} provided that finitely generated subobjects of finitely presented objects are finitely presented.
	\end{itemize}
	In the following proposition we collect some of basic properties of locally finitely presented categories.
	
	\begin{proposition}\label{2.14}
		\begin{itemize}
			\item[(i)]
			Every locally finitely presented abelian category is a Grothendieck category.
			\item[(ii)]
			An object in a locally finitely presented abelian category is finitely
			generated if and only if it is a quotient of some finitely presented object.
			\item[(iii)]
			A locally finitely presented abelian category $\mathcal{A}$ is locally coherent if and only if
			$\fp(\mathcal{A})$ is an abelian category.
			\item[(iv)]
			A full subcategory of a locally finitely presented category is covariantly finite if and only if it is closed under products.
			\item[(v)]
			A full subcategory of a locally finitely presented abelian category that is closed under direct limits and is generated by a set of objects is contravariantly finite.
		\end{itemize}
		\begin{proof}
			We refer the reader to \cite[Pages 1652, 1653 and 1665]{CB} for the statements (i), (ii) and (iv), to \cite[Page 204]{Ro} for (iii), and to \cite[Theorem 3.2]{El} for the last statement.
		\end{proof}
	\end{proposition}
	Let $\mathcal{C}$ be a small additive category. Recall that there is a tensor product bifunctor
	\begin{align}
		\Mod\mathcal{C}\otimes_{\mathcal{C}}& \Mod\mathcal{C}^{op}\longrightarrow \rm Ab. \notag \\
		(F,G)&\longmapsto F\otimes_{\mathcal{C}}G \notag
	\end{align}
	A functor $F\in \Mod\mathcal{C}$ is said to be {\em flat} provided that $F\otimes_{\mathcal{C}}-$ is an exact functor. The subcategory of flat functors is denoted by $\Flatt(\mathcal{C})$. A theorem of Lazard state that an $R$-module is flat if and only if it is a direct limit (over a directed set) of finitely generated free modules. Lazard theorem has been generalized to functors by Oberst and Rohrl \cite{OR}. Indeed a functor $F\in \Mod\mathcal{C}$ is flat if and only if $F$ is a direct limit of representable functors. Some people have called this objects, $\ind$-objects (for more details see \cite[Section 6]{KS}.
	
	The following characterization of locally finitely presented categories is due to Crawley-Boevey \cite{CB}.
	\begin{theorem}$($\cite[Theorem 1.4]{CB}$)$\label{2.15}
		\begin{itemize}
			\item[(a)]
			If $\mathcal{C}$ is an essentially small additive category, then $\Flatt(\mathcal{C})$ is a locally finitely presented category and $\fp(\Flatt(\mathcal{C}))$ consists of direct summands of representable functors. If $\mathcal{C}$ has split idempotents, then the Yoneda functor $h:\mathcal{C}\rightarrow \fp(\Flatt(\mathcal{C}))$ is an equivalence.
			\item[(b)]
			If $\mathcal{B}$ is a locally finitely presented category, then $\fp(\mathcal{B})$ is an
			essentially small additive category with split idempotents and the functor
			\begin{align}
				g:\mathcal{B}&\longrightarrow \Flatt(\fp(\mathcal{B}))\notag \\
				M&\longmapsto \Hom_{\mathcal{B}}(-,M)|_{\fp(\mathcal{B})}\notag
			\end{align}
			is an equivalence.
		\end{itemize}
	\end{theorem}
	
	The following lemma characterize when a functor $F\in \Mod\mathcal{C}$ is a flat functor (equivalently it is an $\ind$-object in $\Mod\mathcal{C}$).
	
	\begin{lemma}\label{2.16}
		Let $\mathcal{C}$ be a small additive category and $F\in \Mod\mathcal{C}$. The following statements are equivalent.
		\begin{itemize}
			\item[(i)]
			$F\in \Ind (\mathcal{C})$.
			\item[(ii)]
			Whenever $X^1\rightarrow X^2$ is a weak cokernel of $X^0\rightarrow X^1$, the sequence $F(X^2)\rightarrow F(X^1)\rightarrow F(X^0)$ is exact.
		\end{itemize}
		\begin{proof}
			See \cite[Page 1644]{CB}.
		\end{proof}
	\end{lemma}
	
	\begin{definition}
		Let $\mathcal{C}$ be a small additive category. An $\ind$-object in $\mathcal{C}$ is an object $F\in \Mod\mathcal{C}$ which is isomorphic to a direct limit of representable functors (i.e. it is a flat functor).
		The full subcategory of $\Mod\mathcal{C}$ consisting of all $\ind$-objects is called the completion (indization in \cite[Section 2.6]{KS}) of $\mathcal{C}$ and is denoted by $\Ind(\mathcal{C})$. So $\Ind (\mathcal{C})= \Flatt(\mathcal{C})$.
	\end{definition}
	Note that the Yoneda functor $\mathcal{Y}:\mathcal{C}\rightarrow \Mod\mathcal{C}$ induces a functor
	\begin{equation}
		\mathbb{I}:\mathcal{C}\rightarrow \Ind (\mathcal{C}). \notag
	\end{equation}
	Obviously $\Ind(\mathcal{C})$ has all filtered limits, and any other functor $\mathcal{C}\rightarrow \mathcal{D}$, where $\mathcal{D}$ has all filtered limits factor uniquely through $\mathbb{I}$.

	%%%%%%%%%%%%%%%%%%%%%%%%%%%%%%%%%%
	
	\section{Completion of $d$-abelian categories}
	
	Let $\mathcal{M}$ be a small $d$-abelian category. In this section we study the completion $\Ind(\mathcal{M})$ of $\mathcal{M}$. In Theorem \ref{3.4} we show that $\Ind(\mathcal{M})$ is "too close" to be a $d$-abelian category.
	
	By Proposition \ref{2.11} we have the following localisation situation of abelian categories.
	\begin{equation}
		\begin{tikzpicture}
			\node (X1) at (-3,0) {$\Eff(\mathcal{M})$};
			\node (X2) at (0,0) {$\Mod\mathcal{M}$};
			\node (X3) at (3,0) {$\frac{\Mod\mathcal{M}}{\Eff(\mathcal{M})}$};
			\node (p) at (-1.5,-0.7) {$p$};
			\node (r) at (1.5,-0.7) {$r$};
			\draw [->,thick] (X1) -- (X2) node [midway,above] {$i$};
			\draw [->,thick] (X2) -- (X3) node [midway,above] {$e$};
			\draw [->,thick] (X2) to [out=225,in=315] (X1) node [midway,left] {};
			\draw [->,thick] (X3) to [out=225,in=315] (X2) node [midway,left] {};
		\end{tikzpicture}
	\end{equation}
	
	\begin{lemma}\label{3.1} Let $\mathcal{M}$ be a small $d$-abelian category. Then
		$\eff(\mathcal{M})=\Eff(\mathcal{M})\cap \modd\mathcal{M}$ and $\Ind(\eff(\mathcal{M}))=\Eff(\mathcal{M})$.
		\begin{proof}
			Let $F\in \eff(\mathcal{M})$. There is a projective presentation
			
			\begin{equation}\label{PP}
				\begin{CD}
					\Hom_\mathcal{M}(-,Y')@ > \Hom_\mathcal{M}(-,f)>> \Hom_\mathcal{M}(-,Y)@>  >>F@ > >>0,
				\end{CD}
			\end{equation}

			induced by an epimorphism $f:Y'\rightarrow Y$. An element $x\in F(X)$ is the image of some morphism $X\rightarrow Y$. Embed $f$ in a $d$-exact sequence and then take $d$-pullback, by \cite[Theorem 3.8]{J} we have the following commutative diagram, where the horizontal morphisms are epimorphic.
			\begin{center}
				\begin{tikzpicture}
					\node (X1) at (-2,2) {$X'$};
					\node (X2) at (0,2) {$X$};
					\node (X3) at (-2,0) {$Y'$};
					\node (X4) at (0,0) {$Y$};
					\draw [->>,thick] (X1) -- (X2) node [midway,above] {$g$};
					\draw [->,thick] (X2) -- (X4) node [midway,above] {};
					\draw [->,thick] (X1) -- (X3) node [midway,left] {};
					\draw [->>,thick] (X3) -- (X4) node [midway,below] {$f$};
				\end{tikzpicture}
			\end{center}
			Obviously $F(g)(x)=0$. Thus $F\in \Eff(\mathcal{M})$.
			
			Now let $F\in \Eff(\mathcal{M})\cap \modd\mathcal{M}$ and choose a projective presentation like \eqref{PP}. Assume that $y$ is the image of $1_Y$ under the map $\Hom_\mathcal{M}(Y,Y)\rightarrow F(Y)$. By the hypothesis there exists an epimorphism $h:X\rightarrow Y$, such that $F(h)(y) =0$. This means that $h$ must factor through $f$. Thus $f$ is an epimorphism.
			
			The last statement follows from the first, because every functor is direct limit of it's finitely presented subfunctors.
		\end{proof}
	\end{lemma}
	
	\begin{proposition}\label{3.2} Let $\mathcal{M}$ be a small $d$-abelian category. Then
		$\frac{\Mod\mathcal{M}}{\Eff(\mathcal{M})}$ is a locally coherent category and $\fp(\frac{\Mod\mathcal{M}}{\Eff(\mathcal{M})})=\frac{\modd\mathcal{M}}{\eff(\mathcal{M})}$.
		\begin{proof}
			It is a well known that $\Mod\mathcal{M}$ is a locally coherent category and $\fp(\Mod\mathcal{M})=\modd\mathcal{M}$. Thus the result follows from Lemma \ref{3.1} and \cite[Proposition A.5]{Kr01}.
		\end{proof}
	\end{proposition}
	
	Let $\mathcal{M}$ be a small $d$-abelian category. By Proposition \ref{2.12}, $\Eff(\mathcal{M})^{_0\bot_d}=\mathcal{L}_d(\mathcal{M})$ is the full subcategory of left $d$-exact functors from $\mathcal{M}$ to the category of abelian groups.
	
	\begin{proposition}\label{3.3}
		Let $\mathcal{M}$ be a small $d$-abelian category. Then we have $\Ind(\mathcal{M})=\mathcal{L}_d(\mathcal{M})$.
		\begin{proof}
			By Yoneda's lemma, every functor $F\in \Mod\mathcal{M}$ is a colimit of representable functors.
			By Lemma 2.16, $F$ is a direct limit of representable functors if and only if whenever $X^1\rightarrow X^2$ is a weak cokernel of $X^0\rightarrow X^1$, the sequence $F(X^2)\rightarrow F(X^1)\rightarrow F(X^0)$ is exact. Thus $\Ind(\mathcal{M})\subseteq\mathcal{L}_d(\mathcal{M})$. Conversely let $F\in \mathcal{L}_d(\mathcal{M})$ and $X^1\rightarrow X^2$ be a weak cokernel of $X^0\rightarrow X^1$. By \cite[Proposition 3.7]{J} there is a right $d$-exact sequence
			\begin{equation}
				X^0\rightarrow X^1\rightarrow X^2\rightarrow X^3\rightarrow \cdots\rightarrow X^d\rightarrow X^{d+1}.\notag
			\end{equation}
			By \cite[Proposition 3.13]{J} we can construct the following commutative diagram
			\begin{center}
				\begin{tikzpicture}
					\node (X1) at (-3,4.5) {$X^0$};
					\node (X2) at (-1,4.5) {$X^1$};
					\node (X3) at (1,4.5) {$\ldots$};
					\node (X4) at (3,4.5) {$X^{d-1}$};
					\node (X5) at (5,4.5) {$X^d$};
					\node (X6) at (7,4.5) {$X^{d+1}$};
					\node (X7) at (-3,3) {$Y_0^1$};
					\node (X8) at (-1,3) {$Y_1^1$};
					\node (X9) at (3,3) {$Y_{d-1}^1$};
					\node (X10) at (5,3) {$Y_d^1$};
					\node (X11) at (7,3) {$Y_{d+1}^1$};
					\node (X12) at (-3,1.5) {$Y_0^2$};
					\node (X13) at (-1,1.5) {$Y_1^2$};
					\node (X14) at (3,1.5) {$Y_{d-1}^2$};
					\node (X15) at (5,1.5) {$Y_d^2$};
					\node (X16) at (7,1.5) {$Y_{d+1}^2$};
					\node (X17) at (-3,0) {$\vdots$};
					\node (X18) at (-1,0) {$\vdots$};
					\node (X19) at (3,0) {$\vdots$};
					\node (X20) at (5,0) {$\vdots$};
					\node (X21) at (7,0) {$\vdots$};
					\node (X22) at (-3,-1.5) {$Y_0^{d-1}$};
					\node (X23) at (-1,-1.5) {$Y_1^{d-1}$};
					\node (X24) at (3,-1.5) {$Y_{d-1}^{d-1}$};
					\node (X25) at (5,-1.5) {$Y_d^{d-1}$};
					\node (X26) at (7,-1.5) {$Y_{d+1}^{d-1}$};
					\node (X27) at (-3,-3) {$Y_0^d$};
					\node (X28) at (-1,-3) {$Y_1^d$};
					\node (X29) at (3,-3) {$Y_{d-1}^d$};
					\node (X30) at (5,-3) {$Y_d^d$};
					\node (X31) at (7,-3) {$Y_{d+1}^d$};
					\node (X32) at (-3,-4.5) {$0$};
					\node (X33) at (-1,-4.5) {$0$};
					\node (X34) at (3,-4.5) {$0$};
					\node (X35) at (5,-4.5) {$0$};
					\node (X36) at (7,-4.5) {$0$};
					\draw [->,thick] (X1) -- (X2) node [midway,above] {};
					\draw [->,thick] (X2) -- (X3) node [midway,above] {};
					\draw [->,thick] (X3) -- (X4) node [midway,above] {};
					\draw [->,thick] (X4) -- (X5) node [midway,right] {};
					\draw [->,thick] (X5) -- (X6) node [midway,right] {};
					\draw [->,thick] (X32) -- (X27) node [midway,right] {};
					\draw [->,thick] (X27) -- (X22) node [midway,right] {};
					\draw [->,thick] (X22) -- (X17) node [midway,right] {};
					\draw [->,thick] (X17) -- (X12) node [midway,right] {};
					\draw [->,thick] (X12) -- (X7) node [midway,above] {};
					\draw [->,thick] (X7) -- (X1) node [midway,above] {};
					\draw [->,thick] (X33) -- (X28) node [midway,above] {};
					\draw [->,thick] (X28) -- (X23) node [midway,right] {};
					\draw [->,thick] (X23) -- (X18) node [midway,right] {};
					\draw [->,thick] (X18) -- (X13) node [midway,right] {};
					\draw [->,thick] (X13) -- (X8) node [midway,right] {};
					\draw [->,thick] (X8) -- (X2) node [midway,right] {};
					\draw [->,thick] (X34) -- (X29) node [midway,right] {};
					\draw [->,thick] (X29) -- (X24) node [midway,above] {};
					\draw [->,thick] (X24) -- (X19) node [midway,above] {};
					\draw [->,thick] (X19) -- (X14) node [midway,above] {};
					\draw [->,thick] (X14) -- (X9) node [midway,right] {};
					\draw [->,thick] (X9) -- (X4) node [midway,right] {};
					\draw [->,thick] (X35) -- (X30) node [midway,right] {};
					\draw [->,thick] (X30) -- (X25) node [midway,right] {};
					\draw [->,thick] (X25) -- (X20) node [midway,right] {};
					\draw [->,thick] (X20) -- (X15) node [midway,right] {};
					\draw [->,thick] (X15) -- (X10) node [midway,above] {};
					\draw [->,thick] (X10) -- (X5) node [midway,above] {};
					\draw [->,thick] (X36) -- (X31) node [midway,right] {};
					\draw [->,thick] (X31) -- (X26) node [midway,right] {};
					\draw [->,thick] (X26) -- (X21) node [midway,right] {};
					\draw [->,thick] (X21) -- (X16) node [midway,right] {};
					\draw [->,thick] (X16) -- (X11) node [midway,above] {};
					\draw [->,thick] (X11) -- (X6) node [midway,right] {};
					\draw [->,thick] (X1) -- (X8) node [midway,above] {};
					\draw [->,thick] (X7) -- (X13) node [midway,right] {};
					\draw [->,thick] (X22) -- (X28) node [midway,right] {};
					\draw [->,thick] (X27) -- (X33) node [midway,right] {};
					\draw [->,thick] (X4) -- (X10) node [midway,right] {};
					\draw [->,thick] (X9) -- (X15) node [midway,right] {};
					\draw [->,thick] (X24) -- (X30) node [midway,above] {};
					\draw [->,thick] (X29) -- (X35) node [midway,right] {};
					\draw [->,thick] (X5) -- (X11) node [midway,right] {};
					\draw [->,thick] (X10) -- (X16) node [midway,right] {};
					\draw [->,thick] (X25) -- (X31) node [midway,right] {};
					\draw [->,thick] (X30) -- (X36) node [midway,right] {};
				\end{tikzpicture}
			\end{center}
			such that for every $k \in \{0, 1, \ldots, d\}$ the diagram
			\begin{center}
				\begin{tikzpicture}
					\node (X1) at (-4,1) {$Y_k^d$};
					\node (X2) at (-2,1) {$Y_k^{d-1}$};
					\node (X3) at (0,1) {$\ldots$};
					\node (X4) at (2,1) {$Y_k^1$};
					\node (X5) at (4,1) {$X^k$};
					\node (X7) at (-4,-1) {$0$};
					\node (X8) at (-2,-1) {$Y_{k+1}^d$};
					\node (X9) at (0,-1) {$\ldots$};
					\node (X10) at (2,-1) {$Y_{k+1}^2$};
					\node (X11) at (4,-1) {$Y_{k+1}^1$};
					\draw [->,thick] (X1) -- (X2) node [midway,above] {};
					\draw [->,thick] (X2) -- (X3) node [midway,above] {};
					\draw [->,thick] (X3) -- (X4) node [midway,above] {};
					\draw [->,thick] (X4) -- (X5) node [midway,above] {};
					\draw [->,thick] (X1) -- (X7) node [midway,above] {};
					\draw [->,thick] (X2) -- (X8) node [midway,left] {};
					\draw [->,thick] (X4) -- (X10) node [midway,left] {};
					\draw [->,thick] (X5) -- (X11) node [midway,left] {};
					\draw [->,thick] (X7) -- (X8) node [midway,above] {};
					\draw [->,thick] (X8) -- (X9) node [midway,above] {};
					\draw [->,thick] (X9) -- (X10) node [midway,above] {};
					\draw [->,thick] (X10) -- (X11) node [midway,above] {};
				\end{tikzpicture}
			\end{center}
			is both $d$-pullback and $d$-pushout diagram. Then its mapping cone is a $d$-exact sequence.
			Applying $F$, we obtain the following commutative diagram
			\begin{center}
				\begin{tikzpicture}
					\node (X0) at (-5,4.5) {$0$};
					\node (X1) at (-3,4.5) {$F(X^{d+1})$};
					\node (X2) at (-1,4.5) {$F(X^d)$};
					\node (X3) at (1,4.5) {$F(X^{d-1})$};
					\node (X4) at (3,4.5) {$\ldots$};
					\node (X5) at (5,4.5) {$F(X^1)$};
					\node (X6) at (7,4.5) {$F(X^0)$};
					\node (X7) at (-3,3) {$F(Y_{d+1}^1)$};
					\node (X8) at (-1,3) {$F(Y_d^1)$};
					\node (X9) at (1,3) {$F(Y_{d-1}^1)$};
					\node (X10) at (5,3) {$F(Y_1^1)$};
					\node (X11) at (7,3) {$F(Y_0^1)$};
					\node (X12) at (-3,1.5) {$F(Y_{d+1}^2)$};
					\node (X13) at (-1,1.5) {$F(Y_d^2)$};
					\node (X14) at (1,1.5) {$F(Y_{d-1}^2)$};
					\node (X15) at (5,1.5) {$F(Y_1^2)$};
					\node (X16) at (7,1.5) {$F(Y_0^2)$};
					\node (X17) at (-3,0) {$\vdots$};
					\node (X18) at (-1,0) {$\vdots$};
					\node (X19) at (1,0) {$\vdots$};
					\node (X20) at (5,0) {$\vdots$};
					\node (X21) at (7,0) {$\vdots$};
					\node (X22) at (-3,-1.5) {$F(Y_{d+1}^{d-1})$};
					\node (X23) at (-1,-1.5) {$F(Y_d^{d-1})$};
					\node (X24) at (1,-1.5) {$f(Y_{d-1}^{d-1})$};
					\node (X25) at (5,-1.5) {$F(Y_1^{d-1})$};
					\node (X26) at (7,-1.5) {$F(Y_0^{d-1})$};
					\node (X27) at (-3,-3) {$F(Y_{d+1}^d)$};
					\node (X28) at (-1,-3) {$F(Y_d^d)$};
					\node (X29) at (1,-3) {$F(Y_{d-1}^d)$};
					\node (X30) at (5,-3) {$F(Y_1^d)$};
					\node (X31) at (7,-3) {$F(Y_0^d)$};
					\node (X32) at (-3,-4.5) {$0$};
					\node (X33) at (-1,-4.5) {$0$};
					\node (X34) at (1,-4.5) {$0$};
					\node (X35) at (5,-4.5) {$0$};
					\node (X36) at (7,-4.5) {$0$};
					\draw [->,thick] (X0) -- (X1) node [midway,right] {};
					\draw [->,thick] (X1) -- (X2) node [midway,above] {};
					\draw [->,thick] (X2) -- (X3) node [midway,above] {};
					\draw [->,thick] (X3) -- (X4) node [midway,above] {};
					\draw [->,thick] (X4) -- (X5) node [midway,right] {};
					\draw [->,thick] (X5) -- (X6) node [midway,right] {};
					\draw [->,thick] (X27) -- (X32) node [midway,right] {};
					\draw [->,thick] (X22) -- (X27) node [midway,right] {};
					\draw [->,thick] (X17) -- (X22) node [midway,right] {};
					\draw [->,thick] (X12) -- (X17) node [midway,right] {};
					\draw [->,thick] (X7) -- (X12) node [midway,above] {};
					\draw [->,thick] (X1) -- (X7) node [midway,above,left] {};
					\draw [->,thick] (X28) -- (X33) node [midway,above] {};
					\draw [->,thick] (X23) -- (X28) node [midway,right] {};
					\draw [->,thick] (X18) -- (X23) node [midway,right] {};
					\draw [->,thick] (X13) -- (X18) node [midway,right] {};
					\draw [->,thick] (X8) -- (X13) node [midway,right] {};
					\draw [->,thick] (X2) -- (X8) node [midway,right] {};
					\draw [->,thick] (X29) -- (X34) node [midway,right] {};
					\draw [->,thick] (X24) -- (X29) node [midway,above] {};
					\draw [->,thick] (X19) -- (X24) node [midway,above] {};
					\draw [->,thick] (X14) -- (X19) node [midway,above] {};
					\draw [->,thick] (X9) -- (X14) node [midway,right] {};
					\draw [->,thick] (X3) -- (X9) node [midway,right] {};
					\draw [->,thick] (X30) -- (X35) node [midway,right] {};
					\draw [->,thick] (X25) -- (X30) node [midway,right] {};
					\draw [->,thick] (X20) -- (X25) node [midway,right] {};
					\draw [->,thick] (X15) -- (X20) node [midway,right] {};
					\draw [->,thick] (X10) -- (X15) node [midway,above] {};
					\draw [->,thick] (X5) -- (X10) node [midway,above] {};
					\draw [->,thick] (X31) -- (X36) node [midway,right] {};
					\draw [->,thick] (X26) -- (X31) node [midway,right] {};
					\draw [->,thick] (X21) -- (X26) node [midway,right] {};
					\draw [->,thick] (X16) -- (X21) node [midway,right] {};
					\draw [->,thick] (X11) -- (X16) node [midway,above] {};
					\draw [->,thick] (X6) -- (X11) node [midway,right]{} ;
					\draw [->,thick] (X7) -- (X2) node [midway,above] {};
					\draw [->,thick] (X12) -- (X8) node [midway,right] {};
					\draw [->,thick] (X27) -- (X23) node [midway,right] {};
					\draw [->,thick] (X32) -- (X28) node [midway,right] {};
					\draw [->,thick] (X8) -- (X3) node [midway,right] {};
					\draw [->,thick] (X13) -- (X9) node [midway,right] {};
					\draw [->,thick] (X28) -- (X24) node [midway,above] {};
					\draw [->,thick] (X33) -- (X29) node [midway,right] {};
					\draw [->,thick] (X10) -- (X6) node [midway,right] {};
					\draw [->,thick] (X15) -- (X11) node [midway,right] {};
					\draw [->,thick] (X30) -- (X26) node [midway,right] {};
					\draw [->,thick] (X35) -- (X31) node [midway,right] {};
				\end{tikzpicture}
			\end{center}
			that satisfies the assumption of the dual of \cite[Lemma 2.9]{EN3}. Thus the top row of the above diagram is an exact sequence and the result follows.
		\end{proof}
	\end{proposition}
	
	In the following theorem we prove some general properties of $\Ind(\mathcal{M})$.
	
	\begin{theorem}\label{3.4}
		Let $\mathcal{M}$ be a small $d$-abelian category. Then the following statements hold.
		\begin{itemize}
			\item[(i)]
			$\Ind(\mathcal{M})$ is a generating-cogenerating functorially finite subcategory of $\frac{\Mod\mathcal{M}}{\Eff(\mathcal{M})}$.
			\item[(ii)]
			For every $F\in \mathcal{L}_1(\mathcal{M})$ there exists
			an exact sequence
			\begin{equation}
				0\rightarrow L_d\rightarrow L_{d-1}\rightarrow \cdots\rightarrow L_1\rightarrow F\rightarrow 0 \notag
			\end{equation}
			in $\mathcal{L}_1(\mathcal{M})$ such that $L_i$'s are in $\mathcal{L}_d(\mathcal{M})$.
			\item[(iii)]
			We have $\Ind(\mathcal{M})^{_1\bot_{d-1}}\subseteq \Ind(\mathcal{M})$.
		\end{itemize}
		\begin{proof}
			(i) Clearly $\Ind(\mathcal{M})$ is a generating subcategory. By Propositions \ref{2.11} and \ref{3.3} we identify $\frac{\Mod\mathcal{M}}{\Eff(\mathcal{M})}$ and $\Ind(\mathcal{M})$ with $\mathcal{L}_1(\mathcal{M})$ and $\mathcal{L}_d(\mathcal{M})$ respectively. By Proposition \ref{2.11}, $\mathcal{L}_1(\mathcal{M})$ is a Grothendieck category and so has enough injective objects. By Lemma \ref{2.12}, every injective object in $\mathcal{L}_1(\mathcal{M})$ is a $d$-exact functor. Thus $\mathcal{L}_d(\mathcal{M})$ is a cogenerating subcategory of $\mathcal{L}_1(\mathcal{M})$. Clearly $\mathcal{L}_d(\mathcal{M})$ is a subcategory of $\mathcal{L}_1(\mathcal{M})$ which is closed under direct limits and products. Therefore by Proposition \ref{2.14}, $\mathcal{L}_d(\mathcal{M})$ is a functorially finite subcategory of locally finitely presented category $\mathcal{L}_1(\mathcal{M})$.
			
			(ii) Let $F\in \mathcal{L}_1(\mathcal{M})$. Choose a right $\mathcal{L}_d(\mathcal{M})$-approximation $L_1\twoheadrightarrow F$ which is an epimorphism by Remark \ref{2.5}. Consider the following short exact sequence of functors in $\Mod\mathcal{M}$.
			\begin{equation}
				0\rightarrow K_2\rightarrow L_1\rightarrow F\rightarrow 0. \notag
			\end{equation}
			Applying this short exact sequence of functors to an arbitrary $d$-exact sequence
			\begin{center}
				$X^0\rightarrow X^1\rightarrow \cdots\rightarrow X^d\rightarrow X^{d+1}$
			\end{center}
			and using the long exact sequence theorem \cite[Theorem 1.3.1]{We} we see that $K_2\in \mathcal{L}_2(\mathcal{M})$. Again choose a right $\mathcal{L}_d(\mathcal{M})$-approximation $L_2\twoheadrightarrow K_2$ and repeating this argument, inductively we obtain an exact sequence
			\begin{equation}
				0\rightarrow L_d\rightarrow L_{d-1}\rightarrow \cdots\rightarrow L_1\rightarrow F\rightarrow 0, \notag
			\end{equation}
			in $\mathcal{L}_1(\mathcal{M})$ with desired properties.
			
			(iii) Let $F\in \Ind(\mathcal{M})^{_1\bot_{d-1}}$. Choose a minimal injective coresolution
			\begin{equation}\label{II}
				0\rightarrow F\rightarrow I^0\rightarrow I^1\rightarrow \cdots\rightarrow I^{d-1}\rightarrow I^d
			\end{equation}
			in $\mathcal{L}_1(\mathcal{M})$. By Proposition \ref{2.10}, $F\in \Ind(\mathcal{M})=\Eff(\mathcal{M})^{_0\bot_{d}}$ if and only if the injective coresolution \eqref{II} is exact in $\Mod\mathcal{M}$. The later follows from the assumption that \eqref{II} remain exact after applying the functor $\Hom_{\Mod\mathcal{M}}(\Hom_\mathcal{M}(-,X),-)$ for every $X\in \mathcal{M}$.
		\end{proof}
	\end{theorem}
	
	By Lemma \ref{2.13}, $\Ind(\mathcal{M})$ is a $d$-rigid subcategory of $\frac{\Mod\mathcal{M}}{\Eff(\mathcal{M})}$ if and only if it is a $d$-rigid subcategory of $\Mod\mathcal{M}$. We conjecture that the following question has positive answer.
	
	\begin{question}\label{3.5}
		Let $\mathcal{M}$ be a small $d$-abelian category. Is $\Ind(\mathcal{M})=\Flatt(\mathcal{M})$ a $d$-rigid subcategory of the abelian category $\Mod\mathcal{M}$?
	\end{question}
	
	\begin{theorem}\label{3.6}
		Assume that Question \ref{3.5} has positive answer. Then
		$\Ind(\mathcal{M})=\mathcal{L}_d(\mathcal{M})$ is a $d$-cluster tilting subcategory of
		the abelian category $\frac{\Mod\mathcal{M}}{\Eff(\mathcal{M})}\simeq\mathcal{L}_1(\mathcal{M})$.
		\begin{proof}
			By Theorem \ref{3.4}, $\Ind(\mathcal{M})$ is a generating-cogenerating functorially finite subcategory of subcategory of
			$\frac{\Mod\mathcal{M}}{\Eff(\mathcal{M})}$ and $\Ind(\mathcal{M})^{_1\bot_{d-1}}\subseteq \Ind(\mathcal{M})$. Because $\Ind(\mathcal{M})$ is $d$-rigid, it remains to show that $^{_1\bot_{d-1}}\Ind(\mathcal{M})\subseteq \Ind(\mathcal{M})$. By Proposition \ref{3.4}, there is an exact sequence
			\begin{equation}
				0\rightarrow L_d\rightarrow L_{d-1}\rightarrow \cdots\rightarrow L_1\rightarrow F\rightarrow 0 \notag
			\end{equation}
			in $\mathcal{L}_1(\mathcal{M})$ such that $L_i$'s are in $\mathcal{L}_d(\mathcal{M})$. Now the result follows from \cite[Proposition 4.4]{Kv2}.
		\end{proof}
	\end{theorem}
	
	\begin{remark}\label{3.7}
		If flat functors and projective functors in $\Mod\mathcal{M}$ are coincide (or equivalently $\mathcal{M}$ be a perfect ring (with several objects)), then $\Ind(\mathcal{M})=\Flatt(\mathcal{M})$ is a $d$-rigid subcategory of the abelian category $\Mod\mathcal{M}$. Note that all $d$-cluster tilting subcategories of finite type (i.e. the number of isomorphism classes of indecomposable objects is finite) of an abelian category $\mathcal{A}$ are perfect.
	\end{remark}
	\section{$d$-cluster tilting subcategories of $\modd A$}
	Let $A$ be a finite-dimensional $K$-algebra, and let $\Mod A$ (res., $\modd A$) be the category of all (res., finitely presented) $A$-modules. If $\mathfrak{M}$ is a $d$-cluster tilting subcategory of $\modd A$, then in "Higher Homological Algebra" we take $\mathfrak{M}$ as a higher-dimensional version of $\modd A$. It is natural to ask the following question: What is the higher-dimensional version of the big module category $\Mod A$? The following proposition shows that if Question \ref{3.5} has positive answer, then the full subcategory $\overrightarrow{\mathfrak{M}}$ of $\Mod A$ consisting of all direct limit of objects from $\mathfrak{M}$ is a good candidate for the higher-dimensional version of the big module category.
	
	\begin{proposition}\label{3.8}
		Assume that Question \ref{3.5} has positive answer.
		Let $\mathfrak{M}$ be a $d$-cluster tilting subcategory of $\modd A$. Then the full subcategory $\overrightarrow{\mathfrak{M}}$ of $\Mod A$ consisting of all direct limit of objects from $\mathfrak{M}$ is the unique $d$-cluster tilting subcategory of $\Mod A$ containing $\mathfrak{M}$.
		\begin{proof}
			$\Mod A$ is a locally finitely presented category with $\fp(\Mod A)=\modd A$. On the other hand we have that $\frac{\Mod\mathfrak{M}}{\Eff(\mathfrak{M})}$ is a locally finitely presented category with
			\begin{equation}
				\fp(\dfrac{\Mod\mathfrak{M}}{\Eff(\mathfrak{M})})=\dfrac{\modd\mathfrak{M}}{\eff (\mathfrak{M})}\simeq \modd A. \notag
			\end{equation}
			Thus by Theorem \ref{2.15}, $\frac{\Mod\mathfrak{M}}{\Eff(\mathfrak{M})}\simeq\Mod A$ and under this equivalence $\Ind(\mathfrak{M})$ maps to $\overrightarrow{\mathfrak{M}}$.
		\end{proof}
	\end{proposition}
	
	\begin{remark}\label{3.9}
		If Question \ref{3.5} has a positive answer, then Proposition \ref{3.8} provides a complete answer to Question 3.6 in \cite{EN2}.
	\end{remark}
	\begin{remark}\label{3.10} In \cite{I3}, Iyama posed the following question, which is still unanswered:
		\begin{itemize}
			\item[(Q1)] Let $A$ be a finite-dimensional $K$-algebra. Does any $d$-cluster tilting subcategory $\mathfrak{M}$ of $\modd A$ with $d \geq 2$ has an additive generator?
		\end{itemize}
		As discussed in \cite[Remark 4.13]{EN2}, the question posed by Iyama in \cite{I3} can be rephrased as follows, which is an equivalent formulation:
		\begin{itemize}
			\item[(Q2)]
			Let $\mathfrak{M}$ be a $d$-cluster tilting subcategory of $\modd A$ with $d \geq 2$. Is $\Add(\mathfrak{M})$ (i.e. the full subcategory of $\Mod A$ consist of all direct summands of direct sums of modules in $\mathfrak{M}$) a $d$-cluster tilting subcategory of $\Mod A$?
		\end{itemize}
		If Question \ref{3.5} has a positive answer, Proposition \ref{3.8} implies that the full subcategory $\overrightarrow{\mathfrak{M}}$ of $\Mod A$ consisting of all direct limits of objects from $\mathfrak{M}$ is a $d$-cluster tilting subcategory of $\Mod A$. Hence, we can rephrase Iyama's question as follows:
		\begin{itemize}
			\item[(Q3)]
			Let $\mathfrak{M}$ be a $d$-cluster tilting subcategory of $\modd A$ with $d \geq 2$. Is any module in $\overrightarrow{\mathfrak{M}}$ a direct sum of finitely presented modules? Equivalently, for every $d$-cluster tilting subcategory $\mathfrak{M}$ of $\modd A$, is $\overrightarrow{\mathfrak{M}}$ a pure semisimple $d$-abelian category (for more details, see \cite{EN2})?
		\end{itemize}
	\end{remark}

Recall that a module $M\in \Mod A$ is said to have a perfect decomposition if there is a decomposition $M=\bigoplus_{k\in \mathcal{K}} X_k$, where $(X_k)_{k\in \mathcal{K}}$ is a locally $T$-nilpotent family of indecomposable modules. This means that for each sequence of non-isomorphisms
\begin{center}
	$X_1 \overset{f_1}{\longrightarrow} X_2 \overset{f_2}{\longrightarrow} X_3\overset{f_4}{\longrightarrow}\cdots$
\end{center}
and each element $x\in X_1$, there
exists $m = m_x$ such that $f_m f_{m-1}\cdots f_1(x) = 0$ (for more details, see \cite{HS}).

\begin{proposition}
	Let $A$ be a finite-dimensional algebra, $\mathfrak{M}$ be a $d$-cluster tilting subcategory of $\modd A$ and $\{{X_i}\}_{i\in \mathcal{I}}$ be a representing set of all indecomposable modules in $\mathfrak{M}$. Put $M:=\bigoplus_{i\in \mathcal{I}} X_i$. Then the following statements are equivalent.
	\begin{itemize}
		\item[(1)] $\mathfrak{M}$ is of finite type.
		\item[(2)] $\Add(\mathfrak{M})=\overrightarrow{\mathfrak{M}}$.
		\item[(3)] $\Add(\mathfrak{M})$ is a covering class in $\Mod A$.
		\item[(4)] $M$ has a perfect decomposition.
	\end{itemize}
\begin{proof}
	It follows from Remark \ref{3.7} and \cite[Corollary 2.3]{S}.
\end{proof}
\end{proposition}

\begin{remark}
	Let $A$ be a finite-dimensional algebra, $\mathfrak{M}$ be a $d$-cluster tilting subcategory of $\modd A$, and $\{{X_i}\}_{i\in \mathcal{I}}$ be a representing set of all indecomposable modules in $\mathfrak{M}$. Motivated by the above result, it is natural to ask whether $M=\bigoplus_{i\in \mathcal{I}}X_i$ is always a perfect decomposition. This question is equivalent to Iyama's question about the finiteness of $\mathfrak{M}$, as shown by the above results. Moreover, by using the equivalent conditions for perfect decomposition \cite[Theorem 1.1]{HS}, this is equivalent to asking whether $\End(N)/J(\End(N))$ is von Neumann regular and idempotents lift modulo $J(\End(N))$ for any $N\in \Add(M)$.
\end{remark}

\begin{remark}
		Let $A$ be a finite-dimensional algebra and $\mathfrak{M}$ be a $d$-cluster tilting subcategory of $\modd A$. We recall that a family of morphisms is called {\it noetherian} if for each sequence $X_1 \overset{f_1}{\rightarrow} X_2 \overset{f_2}{\rightarrow} X_3 \rightarrow \cdots$ of morphisms in the family such that $f_i \circ f_{i-1} \circ \cdots \circ f_1 \neq 0$ for all $i$, there is an integer $n$ such that $f_k$ is an isomorphism for all $k \geq n$. By \cite[Corollary 4.6]{FN}, we know that $\mathfrak{M}$ is of finite type if and only if any family of morphisms between indecomposable modules in $\mathfrak{M}$ is noetherian.
		
		Let
		\begin{center}
			$X_1 \overset{f_1}{\longrightarrow} X_2 \overset{f_2}{\longrightarrow} X_3\overset{f_4}{\longrightarrow}\cdots$
		\end{center}
	be a sequence of morphisms between indecomposable modules in $\mathfrak{M}$. We have the following short exact sequence.
	\begin{center}
		$0\rightarrow \bigoplus_{i\geq 1}X_i\rightarrow \bigoplus_{i\geq 1}X_i \rightarrow \underrightarrow{\Lim}X_i\rightarrow 0.$
	\end{center}
By \cite[Lemma 6.28]{GT}, this short exact sequence splits. Thus $\underrightarrow{\Lim}X_i \in \Add(\mathfrak{M})$. Now, let
		\begin{center}
			$X_1 \overset{f_1}{\longrightarrow} X_2 \overset{f_2}{\longrightarrow} X_3\overset{f_4}{\longrightarrow}\cdots$
		\end{center}
	be a sequence of monomorphisms between indecomposable modules in $\mathfrak{M}$ such that $f_i \circ f_{i-1} \circ \cdots \circ f_1 \neq 0$ for all $i$ and there is an integer $n$ such that $f_k$ is not isomorphism for all $k \geq n$. By \cite{AuL}, we have $\underrightarrow{\Lim}X_i$ is a large module. However, we also have $\underrightarrow{\Lim}X_i \in \Add(\mathfrak{M})$ which is a contradiction. Thus, we conclude that any family of monomorphisms between indecomposable modules in $\mathfrak{M}$ is noetherian. For $\modd A$, this is equivalent to the representation finiteness of $A$ (see \cite[Theorem A]{AuL}). However, it is still unknown whether this implies the finiteness of $\mathfrak{M}$ or not (see \cite[Question 3.4]{DN}).
\end{remark}

	Now we consider the case $d=2$. Recall that a complete cotorsion pair in an abelian category $\mathcal{A}$ is a pair $(\mathcal{X},\mathcal{Y})$ of subcategories of $\mathcal{A}$ that are closed under direct summand and satisfying the following conditions.
\begin{itemize}
	\item[a)] $\Ext^1(X,Y)=0$, for all $X\in \mathcal{X}$ and all $Y\in \mathcal{Y}$.
	\item[b)] For each $M\in \mathcal{A}$, there are exact sequences
	\begin{align*}
		0\rightarrow M\rightarrow Y^{M}\rightarrow X^{M}\rightarrow 0\\
		0\rightarrow Y_{M}\rightarrow X_{M}\rightarrow M\rightarrow 0,
	\end{align*}
such that $Y^{M}, Y_{M}\in \mathcal{Y}$ and $ X^{M}, X_{M}\in \mathcal{X}$.
\end{itemize}

\begin{remark}\label{30.10} Let $A$ be a finite-dimensional algebra. It is well known that a pair $(\mathcal{X},\mathcal{Y})$ is a complete cotorsion pair in an abelian category $\modd A$ if and only if it satisfies the following conditions (see \cite[Section 5.2]{GT}).
	\begin{itemize}
		\item[a)] $\mathcal{X}=^{\perp_1}\mathcal{Y}$ and $\mathcal{Y}=\mathcal{X}^{\perp_1}$.
		\item[b)] In the above short exact sequences $M\rightarrow Y^M$ is a left $\mathcal{Y}$-approximation and $X_{M}\rightarrow M$ is a right $\mathcal{X}$-approximation.
		\begin{align*}
		\end{align*}
	\end{itemize}
\end{remark}

By the above remark, we conclude the following obvious observation.

\begin{proposition}\label{31.10}
	Let $A$ be a finite-dimensional algebra. A subcategory $\mathfrak{M}$ of $\modd A$ is a $2$-cluster tilting subcategory if and only if $(\mathfrak{M},\mathfrak{M})$ is a complete cotorsion pair.
\end{proposition}
\begin{proposition}\label{32.10}
	Let $A$ be a finite-dimensional algebra and $\mathfrak{M}$ be a $2$-cluster tilting subcategory of $\modd A$. Then $(\Add(\mathcal{M}),\overrightarrow{\mathfrak{M}})$ is a complete cotorsion pair in $\Mod A$.
	\begin{proof}
		By \cite[Theorem 6.11]{GT}, $(^{\perp}(\mathfrak{M}^{\perp}),\mathfrak{M}^{\perp})$ is a complete cotorsion pair in $\Mod A$. By the proof of Theorem \ref{3.4} $(iii)$ we know that $\mathfrak{M}^{\perp}=\overrightarrow{\mathfrak{M}}$, and by \cite[Corollary 6.14]{GT} and \cite[Lemma 7.2]{HHT}, $^{\perp}(\mathfrak{M}^{\perp})=\Add(\mathfrak{M})$.
	\end{proof}
\end{proposition}

Now, we can prove the following result, which shows that when $d=2$, Question \ref{3.5} is equivalent to Iyama's Question.

\begin{theorem}\label{4.7}
	Let $A$ be a finite-dimensional algebra and $\mathfrak{M}$ be a $2$-cluster tilting subcategory of $\modd A$. The following are equivalent.
	\begin{itemize}
		\item[(1)] $\mathfrak{M}$ is of finite type.
		\item[(2)] $\Add(\mathfrak{M})=\overrightarrow{\mathfrak{M}}$.
		\item[(3)] $\overrightarrow{\mathfrak{M}}$ is a $2$-cluster tilting subcategory of $\Mod A$.
		\item[(4)] $\Add(\mathfrak{M})$ is a $2$-cluster tilting subcategory of $\Mod A$.
	\end{itemize}
\begin{proof}
	 $(2)$ is equivalent to $(3)$ and also to $(4)$ by Propositions \ref{31.10} and \ref{32.10}. If $\mathfrak{M}$ is of finite type, then by \cite[Theorem 3.3]{EN2}, $\Add(\mathfrak{M})$ is a $2$-cluster tilting subcategory of $\Mod A$. If $\Add(\mathfrak{M})$ is a $2$-cluster tilting subcategory of $\Mod A$, then by \cite[Theorem 4.12]{EN2}, $\mathfrak{M}$ is of finite type.
\end{proof}
\end{theorem}

	\section*{acknowledgements}
	The research of the first author was in part supported by a grant from IPM (No. 1400180047). The research of the second author was in part supported by a grant from IPM (No. 14020416).
	The work of the authors is based upon
	research funded by Iran National Science Foundation (INSF) under project No. 4001480.

\end{document}